\documentstyle[11pt]{article}

\input{psfig}
\begin{document}
\centerline {\Large \bf Dynamics of Certain Smooth
One-dimensional Mappings}

\large

\vskip5pt
\centerline {\Large \bf I. The $C^{1+\alpha}$-Denjoy-Koebe distortion lemma}

\vskip10pt
\centerline{Yunping Jiang }
\centerline{Institute for Mathematical Sciences, SUNY at Stony Brook}
\centerline{Stony Brook, L.I., NY 11794}
\vskip5pt
\centerline{June 24, 1990}
\vskip25pt
\centerline{ {\Large \bf Abstract}}
\vskip5pt
We prove a technical lemma, the $C^{1+\alpha }$-Denjoy-Koebe distortion lemma,
estimating the distortion of a long composition
of a $C^{1+\alpha }$ one-dimensional mapping $f:M\mapsto M$
with finitely many, non-recurrent,
power law critical points. The proof of this lemma combines the ideas of
the distortion lemmas of
Denjoy and Koebe.

\vskip25pt
\centerline{{\Large \bf Contents}}
\vskip5pt

\noindent {\large \S 1 Introduction.}

\vskip5pt
\noindent {\large  \S 2 A Very Good Mapping.}

\vskip3pt
\S 2.1 a power law critical
point.

\vskip3pt
\S 2.2 The new differentiable structure
associated with a semi-good

mapping.

\vskip3pt
\S 2.3 The definition of a very good mapping.

\vskip5pt
\noindent {\large \S 3
The Distortion Of A Long Composition Of A
Very Good Mapping.}

\vskip3pt
\S 3.1 The naive distortion lemma.

\vskip3pt
\S 3.2 The proof of $C^{1+\alpha }$-Denjoy-Koebe distortion
lemma.

\vskip5pt
\S 3.3 A larger class of one-dimensional mappings.

\pagebreak

\noindent {\Large {\bf \S 1 Introduction}}
\vskip5pt

There are two techniques in studying the distortion of a long
composition
of a one-dimensional smooth mapping.

\vskip5pt
\noindent {\bf ``Denjoy Principle'':} One technique
goes back to Denjoy. Many people have contributed to this technique
[D], [S], [N1], [N2], [N3], [H], [M], etc.. We
call one of the formulations of this technique \underline
{the naive} \underline{distortion
lemma} because its proof is straightforward -- any one, who has
been trained in Calculus, will understand the proof ([J1],
p25--26, or Lemma 3 in this paper).
The naive distortion lemma is one of the key lemmas in studying
a long composition of a one-dimensional $C^{1+\alpha }$-endomorphism.

\vskip5pt
\noindent {\bf ``Koebe principle'':} The second technique was found in
recent years in studying a long composition of a mapping with
critical points from a one-dimensional manifold to itself.  Many people
formulated this principle in different ways, [GS], [Su1], [Su2], [WS],
[Sw], etc.. We call
one version  \underline {the $C^{3}$-Koebe
distortion
lemma} (see also [J1, p26--27] for a complete proof).  I learned
this
from Sullivan, who invented the name ``Koebe principle''
in analogy with the Koebe lemma in one variable
complex analytic functions.
We consider the nonlinearity of a
$C^{2}$-function $f$ on an interval $I$ as a one-form
\[ n(f)= \frac{f''}{f'}dx .  \]
If the nonlinearity of the function $f$ is integrable on $I$, then the
distortion $|f'(x)/f'(y)|$ of $f$ at any pair $x$ and $y$ in $I$ is
bounded.
The problem is that the nonlinearity of $f$ may be non-integrable if $f$
has a critical point. The $C^{3}$-Koebe distortion
lemma estimates the nonlinearity of a one-dimensional $C^{3}$-mapping
$f$ with nonnegative Schwarzian derivative. This property,
nonnegative Schwarzian derivative, is preserved under
composition, which makes the $C^{3}$-Koebe distortion lemma a very
useful technique in studying a long composition of a one-dimensional
$C^{3}$-mapping with nonpositive Schwarzian derivative (its inverse branches
have nonnegative Schwarzian derivatives). However, the assumption of
nonpositive Schwarzian derivative is a very strong one.

\vskip5pt
{\bf What we would like to say in this paper.}	We
prove a technical lemma, \underline {the $C^{1+\alpha }$-Denjoy-Koebe distortion
lemma}, estimating the distortion of a long
composition
of a one-dimensional $C^{1+\alpha}$-mapping
with finitely many non-recurrent critical points of certain types.
The formulation and the proof of this lemma combine the ideas
of the distortion lemmas of Denjoy and
Koebe.

\vskip5pt
Suppose $M$ is an oriented connected one-dimensional
$C^{2}$-Riemannian manifold with Riemannian metric $dx^{2}$ and associated
length element $dx$. Suppose $f:M\mapsto M$ is a continuous mapping.
A \underline {critical point} $c$ of $f$ is a point in $M$ such that
either $f$ is not differentiable at this point or $f$ is
differentiable at this point with zero derivative. We always assume
that $f$ is $C^{1}$ at a noncritical point $p$, namely there
is a neighborhood $U_{p}$ of $p$ such that the restriction of $f$ to
$U_{p}$ is differentiable
and the derivative $(f|U_{p})'$ is continuous.
We call the
image of a critical point under $f$ a \underline{critical value} of $f$.
We say a critical point $c$ of $f$ is a \underline{power} \underline{law} critical point if it
is an isolated critical point and there is a number $\gamma\geq 1$ such
that the limits of ratio, $f'(x)/|x-c|^{\gamma -1}$, exist and equal
nonzero numbers as $x$ goes to $c$ from below and from above. We call
the number $\gamma$ the \underline{exponent} of $f$ at the power law critical point
$c$. We will assume that $f:M\mapsto M$ is a $C^{1}$-mapping for we are only
interested in a smooth critical point of $f$.

\vskip5pt
For a $C^{1}$-mapping $f: M\mapsto M$ with only power law critical points
such that
the set of critical points and the set of critical values of $f$
are disjoint, we
define
a \underline{new differentiable structure} on the underlying space
$M$.
This new differentiable structure associated with the mapping $f$
has the local parameter
$\int dx/|x|^{\tau}$, where $\tau =1-1/\gamma$, on a neighborhood of a
critical value of $f$
if the corresponding critical point has the exponent $\gamma$.
On a neighborhood of any other point, the new
differentiable structure has the local parameter $\int \rho (x) dx$
where $\rho (x)$ is a positive $C^{2}$-function.
{\it With respect to the new differentiable
structure, the left and the right derivatives of $f$ at
any critical point exist and equal
nonzero numbers} (see Figure 1). We call the original
differentiable structure the old one.

\begin{figure}
\centerline{\psfig{figure=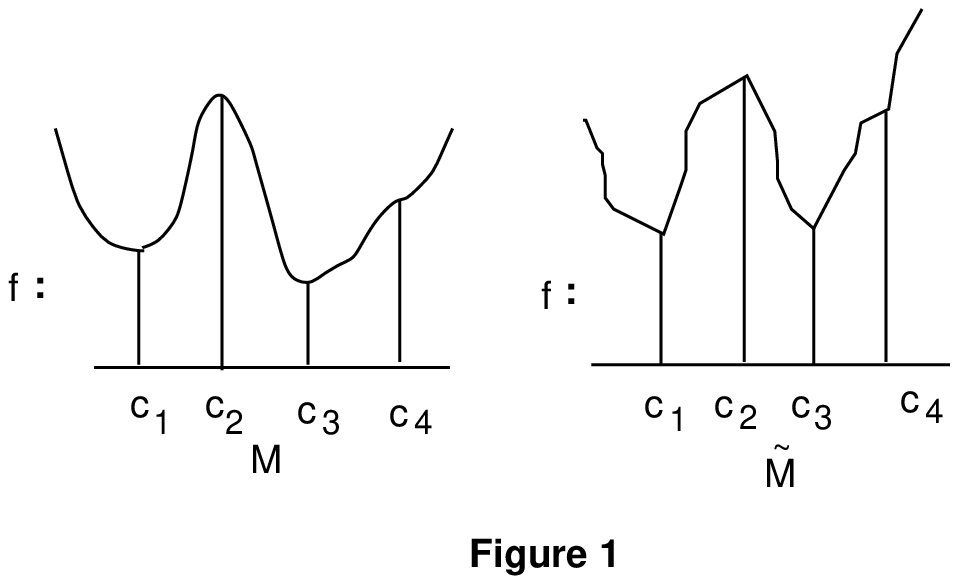}}
\end{figure}

\vskip5pt
We use the oriented connected one-dimensional smooth manifolds
$M$ and $\tilde{M}$,
which
are the same topological space but with the old and the new differentiable
structures , respectively,
to study the dynamics of the mapping $f:M\mapsto M$:
the distortions of a long composition of a one-dimensional
$C^{1+\alpha}$-mapping $f
: M\mapsto M$
with only finitely many, non-recurrent, power law critical points
has an estimate like that in the naive distortion lemma and that in
the $C^{3}$-Koebe distortion lemma.

\vskip20pt

{\bf Acknowledgement.} The preparation of this manuscript was completed
in the Graduate Center of CUNY and the IMS in SUNY at Stony Brook. It is
pleasure for me to thank D. Sullivan for many useful
discussions and G. Swiatek and E. Cawley for reading the
manuscript. I want to thank J. Milnor for reading and correcting the
manuscript and for his very helpful remarks, suggestions and criticisms of
this and other my manuscripts.

\vskip20pt
\centerline {\Large {\bf \S 2
A Very Good Mappings}}

\vskip7pt
Suppose $M$ is an oriented connected one-dimensional
$C^{2}$-Riemannian manifold with Riemannian metric $dx^{2}$ and associated
length element $dx$. Suppose
$f: M\mapsto M$ is a continuous mapping.
A critical point $c$ of $f$ is a point in $M$ such that

\vskip5pt
$(a)$ $f$ is not differentiable at this point or

\vskip5pt
$(b)$ $f$ is
differentiable at this point but the derivative of $f$ at this point is
zero.

\vskip7pt
We always assume
that $f$ is $C^{1}$ at a noncritical point $p$, namely there
is a neighborhood $U_{p}$ of $p$ such that the restriction of $f$ to
$U_{p}$ is differentiable and
the derivative $(f|U_{p})'$ is
continuous.
We call the
image of a critical point under $f$ a critical value of $f$.

\vskip7pt
\noindent {\large \bf  \S 2.1  A power law critical
point.}

\vskip7pt
We give a definition of a power law critical point for the
one-dimensional mapping $f: M\mapsto M$ as follows.

\vskip7pt
{\sc Definition 1.}
{\em Suppose $c$ is an isolated critical point of $f$ and
suppose there are $\gamma^{-}$,
$\gamma^{+}\geq 1$ such that
\[ \lim_{x\mapsto c-}f'(x)/|x-c|^{\gamma^{-} -1}=A
\hskip5pt and \hskip5pt
\lim_{x\mapsto c+}f'(x)/|x-c|^{\gamma^{+} -1}=B\]
exist and equal nonzero numbers.
Then we say that $c$ is a power law critical point with
the left and right exponents $\gamma^{-}$ and $\gamma^{+}$.}

\vskip10pt
\centerline{\psfig{figure=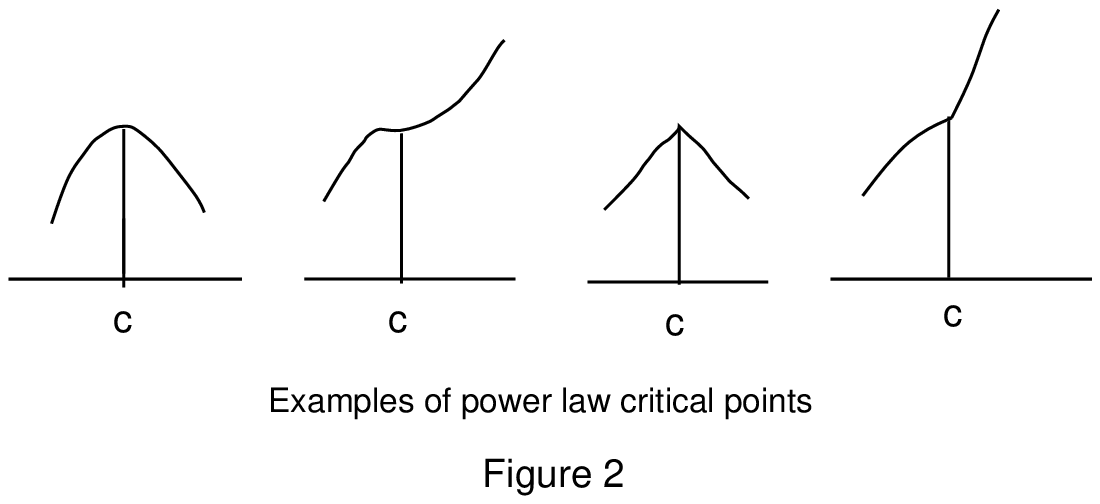}}

\vskip5pt
The following is essentially proved in [J4] (see [J5], too).

\vskip5pt
{\sc Preliminary Lemma.}
{\em Suppose $f: M\mapsto M$ is a continuous mapping
and $c$
is a power law critical point with the left and right exponents
$\gamma^{-}$ and $\gamma^{+}$. Then
there is a continuous mapping $\tilde{f}: M\mapsto M$
and a real number $\sigma \neq 0$ such that

\vskip3pt
$(a)$ the mapping $\tilde{f}$ either has the form

\[ \tilde{f}=\left\{
\begin{array}{ll}
  -\sigma |x-c|^{\gamma^{-}}+f(c) & x\leq c, \\
  |x-c|^{\gamma^{+}}+f(c) & x\geq c
\end{array}
\right.
\hskip5pt or \hskip5pt
\tilde{f}=\left\{
\begin{array}{ll}
  \sigma |x-c|^{\gamma^{-}}+f(c) & x\leq c, \\
  -|x-c|^{\gamma^{+}}+f(c) & x\geq c
\end{array}
\right.
\]
where $x$ is in a small
neighborhood of $c$,

\vskip3pt
$(b)$ the mapping $f$ is semi-conjugate to the mapping $\tilde{f}$.
This means that there is a
monotone and continuous mapping $h$ from $M$ onto $M$
such that
\[ h\circ f=\tilde{f}\circ h\]
and $h$ is differentiable at $c$ with $h'(c)>0$.

\vskip3pt
Moreover,

$(i)$ if $f$ is $C^{1+\alpha }$ on $x\leq c$ and on $x\geq c$
for some
$0< \alpha \leq 1$, and $r_{-}(x)=f'(x)/|x-c|^{\gamma^{-} -1}$, $x\leq c$,
 and
$r_{+}(x)=f'(x)/|x-c|^{\gamma^{+} -1}$, $x\geq c$,  are $C^{\beta}$ for some
$0< \beta \leq 1$, where $x$ is in a small neighborhood of $c$, then the
mapping $h$ can be an orientation-preserving $C^{1}$-diffeomorphism.

\vskip3pt
$(ii)$
The mapping $h$ can be
an orientation-preserving $C^{1,1}$ or $C^{2}$-diffeomorphism if and only if
$f$ is $C^{1,1}$ or $C^{2}$ on $x\leq c$ and $x\geq c$, and
$r_{-}(x)$ and
$r_{+}(x)$ are Lipschitz or $C^{1}$, where $x$ is in a small
neighborhood of $c$. }

\vskip5pt
{\bf Remark.}
 For a power law critical point of $f$, the left and right exponents
are $C^{1}$-invariants. By this we mean that they are the same numbers for $f$ and for
$h\circ f\circ
h^{-1}$ whenever $h$ is an orientation-preserving $C^{1}$-diffeomorphism.
When the left and right exponents are the same, we then have an
important $C^{1}$-invariant
\[ \sigma =\lim_{x\mapsto c-}\frac{f'(x)}{f'(-x+2c)}\]
which we call the asymmetry of $f$ at $c$. The number $\sigma$ in
Preliminary Lemma is the asymmetry. In the paper [J1], we showed
that the asymmetry is an independent $C^{1}$-invariant.

\vskip5pt
\noindent {\large \bf \S 2.2  The new differentiable structure
associated with a semi-good mapping.}

\vskip5pt
Although the results in the rest of the paper hold for a mapping $f$ with both
smooth and non-smooth critical points, but we are only interested in a
smooth critical point of $f$. Henceforth we will assume that $f:M\mapsto M$ is a
$C^{1}$-mapping. Moreover we will assume that the left and right
exponents of $f$
at a power law
critical point are the same.

\vskip5pt
{\sc Definition 2.}
{\em  We say $f$ is a
semi-good mapping if

\vskip5pt
$(I)$ the mapping $f$ has only finitely many power law critical points,

\vskip3pt
$(II)$ the set of
critical points and the set of
critical values
of $f$ are disjoint, and if

\vskip3pt
$(III)$
the exponents of $f$ at two critical points are the same
whenever the images of these two points under $f$
are the same.}

\vskip5pt
Suppose $f: M\mapsto
M$ is a semi-good mapping. Let $CP=\{ c_{1},
\cdots , c_{d}\}$ be the set of critical points of $f$ and
$\Gamma =\{ \gamma_{1}, \cdots \gamma_{d}\}$ be the set of corresponding
exponents.
We define a new differentiable structure associated with
$f$ as follows.

\vskip5pt
Suppose $\Phi=\{ (w_{j}, W_{j}) \}_{j\in \Lambda}$ is a
$C^{2}$-atlas of $M$, this means that $\{ W_{j}\}_{j\in \Lambda}$
is a cover of open sets of $M$ and $\{ w_{j}:W_{j}\mapsto {\bf
R}^{1}\}_{j\in \Lambda}$ is a
set of homeomorphisms such that every $w_{jk}=w_{j}\circ w^{-1}_{k}$
is a $C^{2}$-function whenever $W_{j}$ and $W_{k}$
are overlap.
Suppose every
critical value $v_{i}=f(c_{i})$ is in one and only one chart
$(w_{i}, W_{i})$ and $w_{i}$ maps the critical value $v_{i}$ to $0$.
For every critical value $v_{i}=f(c_{i})$,
we use $k_{i}(x)$ to denote the homeomorphism $\int_{0}^{x}
dx/|x|^{\tau_{i}}: {\bf
R}^{1}\mapsto {\bf R}^{1}$ where $\tau_{i}=1-1/\gamma_{i}$.
Let
$\tilde{w}_{j}=k_{i}\circ w_{i}$ if $W_{j}=W_{i}$ contains a critical
value $v_{i}=f(c_{i})$ and $\tilde{w}_{j}=w_{j}$ if $W_{j}$ does not contain
any critical values.
The set $\tilde {\Phi }=\{ (\tilde{w}_{j}, W_{j})\}$ is another
$C^{2}$-atlas of $M$.
We call the maximal $C^{2}$-atlas of $M$ which contains the set
$\tilde{ \Phi}= \{
(\tilde{w}_{j}, W_{j})\}$ the new differentiable structure associated
with $f$ on $M$.
We denote the topological space $M$ equipped with
this new differentiable structure as a differentiable manifold
$\tilde{M}$.	

\vskip5pt
It is often convenient to
think the new differentiable structure associated with $f$
as a singular metric $\rho (x)dx$ with respect to $dx$ and
the mapping $h=\int \rho (x)dx: M\mapsto M$
as the corresponding
change of coordinate on $M$.
The mapping $\tilde{f}= h\circ f\circ h^{-1}: M\mapsto M$
is the representation of the mapping $f: \tilde{M}\mapsto
\tilde{M}$.

\vskip5pt
{\sc Lemma 1.} {\em Suppose $f: M\mapsto M$ is a semi-good
mapping and $CP$ is the set of critical points of $f$.
Then the mapping $f: \tilde{M}\mapsto \tilde{M}$ is a continuous
mapping
and at every point $c_{i}\in CP$, the left and right
derivatives of $f:\tilde{M}\mapsto \tilde{M}$ exist and equal nonzero
numbers.}

\vskip5pt

{\it Proof.} The proof of this lemma is easy.
The reader may
do it as an exercise or refer to the proof in [J1, p21].

\vskip5pt
\noindent {\large \bf \S 2.3 The definition of a very good mapping.}

\vskip5pt
We define a \underline {very good
mapping.}
Before to give the
definition of a very good mapping, we define
the term $C^{1+\alpha }$ for a real number $0< \alpha \leq 1$ and a semi-good
mapping.

\vskip5pt
Suppose $f:M\mapsto M$ is a semi-good mapping.
Let $CP$ be the set of critical points of $f$. Suppose
$\eta_{0}$ is the set of the closures of the intervals of the
complement of $CP$.  We say a homeomorphism $g:I\mapsto J$ is a
$C^{1+\alpha}$-embedding for some $0< \alpha \leq 1$ if $g$ and $g^{-1}$ are both differentiable
with $\alpha$-H\"older continuous derivatives.

\vskip5pt
{\sc Definition 3.}
{\em we say $f$ is $C^{1+\alpha }$ for some $0< \alpha \leq 1$ if

\vskip3pt
$(1)$ the restriction of $f$ to every interval in $\eta_{0}$
is differentiable with $\alpha$-H\"older continuous derivative,

\vskip3pt
$(2)$ for every critical point $c_{i}$, there is a
neighborhood $U_{i}$ of $c_{i}$ such that the restrictions of
$f:\tilde{M}\mapsto \tilde{M}$ to the intersection of $U_{i}$ and $\{ x\leq
c_{i}\} $ and the intersection of $U_{i}$ and $\{ x\geq c_{i}\}$ are
$C^{1+\alpha }$-embeddings.}

\vskip5pt
We will assume that $U_{i}$ is a closed interval for every $i=1$,
$\cdots$, $d$.	Suppose ${\cal U}$ be the union $\cup_{i=1}^{d}U_{i}$
and ${\cal V}$ be the closure of the complement of ${\cal U}$ in $M$.

\vskip5pt
{\sc Definition 4.} {\em A $C^{1}$-mapping $f: M\mapsto M$
is a very good $C^{1+\alpha}$-mapping (or a very good mapping) for some
$0< \alpha \leq 1$ if it is a semi-good mapping and satisfies

\vskip3pt
$(IV)$
$f$ is $C^{1+\alpha }$,

\vskip3pt
$(V)$ the set $CP$ of critical points and the closure of the
post-critical orbits $\cup_{n=1}^{\infty}f^{\circ n}(CP)$ are
disjoint and

\vskip3pt
$(VI)$ there are two constants $K>0$ and $\nu >1$ such that
for any ${\cal O}_{x, n}=\{ x$, $f(x)$, $\cdots $, $f^{\circ
(n-1)}(x)\}$
with ${\cal O}_{x,n}\cap {\cal U}=\emptyset$,
$|(f^{\circ
k})'(x)|\geq K\nu^{k}$ for any $1\leq k\leq n$.}

\vskip5pt
The space of good mappings is a quite large one,
for example, it contains all $C^{3}$ semi-good mappings with
nonpositive Schwarzian derivative and finitely many non-recurrent
critical points (see, for example, [Mi], [MS] and [J2]).

\vskip10pt
\centerline {\Large {\bf  \S 3
The Distortion Of A Long Composition Of}}

\centerline{\Large {\bf A Very Good Mapping}}

\vskip5pt
Suppose $f: M\mapsto M$ is a very good $C^{1+\alpha}$-mapping for some
$0< \alpha \leq 1$. We always
assume that ${\cal U}$, the union of all $U_{i}$ in Definition 3, is
disjoint with the closure of the post-critical orbits
$\cup_{n=1}^{\infty}f^{\circ n}(CP)$.
We use $U_{i-}$ to
denote the subset consisting of all points $x$ in $U_{i}$ with $x\leq c_{i}$
and use
$U_{i+}$ to
denote the subset consisting of all points $x$ in $U_{i}$ with $x\geq c_{i}$.
Let ${\cal W}$ be the collection of all $U_{i-}$ and $U_{i+}$. Remember
that ${\cal V}$ is the closure of the complement of ${\cal U}$ in $M$.
We say a sequence ${\cal
I}=\{ I_{j}\}_{j=0}^{n}$ of intervals of $M$ is \underline {suitable} if

\vskip3pt
$(i)$ $I_{j}$ is
the image of $I_{j+1}$ under $f$ for $j=0$, $\cdots$ $n-1$ and

\vskip3pt
$(ii)$ either $I_{j}$ is in ${\cal V}$
or $I_{j}$ is in some interval in ${\cal W}$ for every $j=0$, $\cdots
$, $n$.

\vskip5pt
For a suitable
sequence ${\cal I}=\{ I_{j}\}_{j=0}^{n}$ of intervals of $M$,
we use $g_{j}$ to denote the inverse of the restriction of
$f^{\circ j}$ to $I_{j}$. For a pair of points $x$ and $y$ in
$I_{0}$, we use $x_{j}$ and $y_{j}$ to denote the images of $x$ and $y$
under $g_{j}$ and call the ratio
$|g_{n}'(x)|/|g_{n}'(y)|$
the distortion of $f$
at $x$ and $y$ along ${\cal I}$.  We use $D_{xy}$ to denote
the distance between $\{ x, y\}$ and post-critical orbits
$\cup_{n=1}^{\infty } f^{\circ n}(CP)$.

\vskip5pt
The main result of this paper is the
following:

\vskip5pt
{\sc Lemma 2} (the $C^{1+\alpha }$-Denjoy-Koebe distortion lemma).
{\em Suppose $f:M\mapsto M$ is a very good $C^{1+\alpha}$-mapping for
some $0< \alpha \leq 1$. There are two positive constants $A$ and $B$
such that for any
suitable sequence ${\cal I}=\{ I_{j}\}_{j=0}^{n}$ of intervals of $M$
and any pair $x$ and $y$ in $I_{0}$,
the distortion of $f$ at $x$ and $y$ along ${\cal I}$
satisfies
 \[
\frac{|g_{n}'(x)|}{|g_{n}'(y)|}\leq \exp
\Big( A\sum_{i=0}^{n}|x_{i}-y_{i}|^{\alpha}
+ \frac{B|x-y|}{D_{xy}} \Big).\]}

\vskip5pt
\noindent {\large \bf \S 3.1. The naive distortion lemma.}

\vskip5pt
Before to prove Lemma 2, we state the naive distortion lemma.
Suppose $g:U\mapsto M$ is a
$C^{1+\alpha }$-mapping
for some $0< \alpha \leq 1$ where $U$ is an interval of
$M$. Let $K$ be
the $\alpha$-H\"older constant of the derivative of $g$, this means that
$K$ is the smallest positive constant such that
\[ |g'(x)-g'(y)|\leq K|x-y|^{\alpha }\]
for all $x$ and $y$ in $U$.
Suppose $\{
I_{j}\}_{j=1}^{n}$ is a sequence of intervals of $U$ and
$x_{i}$
and $y_{i}$ are two points in $I_{j}$ for $1\leq j\leq n$.
We also call the product of ratios
$\prod_{j=1}^{n}|g'(x_{j})|/ |g'(y_{j})|$ the distortion of
$g$ at
$\{ x_{j}\}_{j=1}^{n}$ and
$\{ y_{j}\}_{j=1}^{n}$.
Let $\beta$ be the minimum of $|g'|$ on $\cup_{j=0}^{n} I_{j}$.

\vskip5pt
{\sc Lemma 3} (the naive distortion lemma).
{\em
The distortion of $g$
at $\{ x_{j}\}_{j=1}^{n}$ and
$\{ y_{j}\}_{j=1}^{n}$ satisfies
\[\prod_{j=1}^{n}\frac{|g'(x_{j})|}{
|g'(y_{j})|}\leq
exp \Big (\frac{K}{\beta}\sum_{j=0}^{n} |x_{j}-y_{j}|^{\alpha }
\Big).\]}

\vskip5pt

{\it Proof.}
Take the function $\log x$ at
$\prod_{j=1}^{n}|g'(x_{j})|/
|g'(y_{j})|$, we have

\[\log \Big( \prod_{j=1}^{n}\frac{|g'(x_{j})|}
{|g'(y_{j})|} \Big) =
\sum_{j=1}^{n} \Big( \log |g'(x_{j})|-\log
|g'(y_{j})| \Big) .\]
Because
$\log x$ is Lipschitz continuous with the
Lipschitz constant $1/\beta$ on the interval $[\beta, +\infty )$ and
the $\alpha$-H\"older constant of $g'$ on $U$ is $K$, we
have that
\[ |\sum_{j=0}^{n} \Big( \log |g'(x_{j})|-\log
|g'(y_{j})| \Big) |\leq
\frac{1}{\beta}\sum_{j=0}^{n} |g'(x_{j})-g'(y_{j})|\]
which is bounded above by
$(K/\beta)\sum_{j=0}^{n} |x_{j}- y_{j}|^{\alpha }$.

\vskip5pt
\noindent {\large \bf \S 3.2 The proof of $C^{1+\alpha }$-Denjoy-Koebe distortion
lemma.}

\vskip5pt
We call ${\cal U}$, the union of $U_{i}$ for $i=1$, $\cdots $, $d$,
the critical set and ${\cal V}$, the closure of the complement of
${\cal U}$ in $M$, the
noncritical set (see Figure 3). Let $\eta_{0}$ be the set of the
closures of the intervals of the complement of the set $CP$ of critical
points of $f$ in $M$.
Let $\tilde{f}=h\circ f\circ h^{-1}:M\mapsto M$ be the representation of
$f:\tilde{M}\mapsto \tilde{M}$, where $h$ is the corresponding change of
coordinate. Remember that $U_{i-}=U_{i}\cap \{ x: x\leq c_{i}\}$ and
$U_{i+}=U_{i}\cap \{x:x\geq c_{i}\}$.

\vskip5pt

\centerline{\psfig{figure=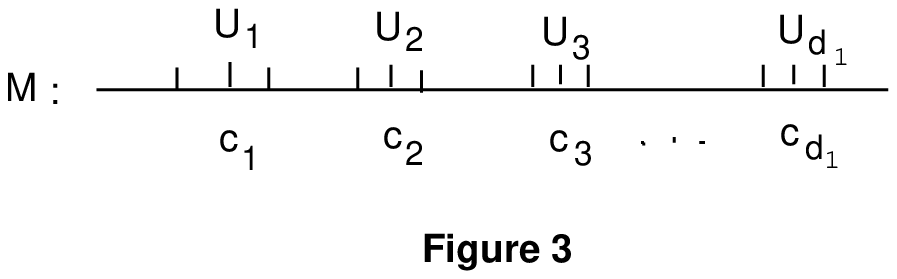}}

\vskip5pt
Let $K_{1}>0$ be the maximum of the $\alpha$-H\"older constants
of the derivatives of the restrictions of $f$ to the intervals
in $\eta_{0}$ and
$\beta_{1}>0$ be the minimum of the absolute value of the restriction of
the derivative $f'$ of $f$ to ${\cal V}$.

\vskip5pt
The restrictions of $\tilde{f}$ to the sets $U_{i-}$ and $U_{i+}$ are
$C^{1+\alpha }$-embeddings for $i=1$, $\cdots $, $d$.
Let
$K_{2}>0$ be the maximum of the $\alpha $-H\"older constants of the
derivatives of these restrictions and $\beta_{2}>0$ be the minimum
of the absolute value of the derivatives of these restrictions.

\vskip5pt
The restrictions of $h$ to the intervals of ${\cal U}$
are $C^{1,1}$.	Let $K_{3}>0$ be the maximum of Lipschitz constants
of the derivatives of these restrictions
and $\beta_{3}>0$ be the minimum of the absolute value of the
derivatives of these restrictions.

\vskip5pt
The distortion of $f$ along ${\cal I}$ at $x$ and $y$ satisfies
\[ \frac{|g_{n}'(x)|}{|g_{n}'(y)|}
=\frac{|(f^{\circ n})'(y_{n})|}{|(f^{\circ n})'(x_{n})|}.\]
By the chain rule,
the ratio $|(f^{\circ n})'(y_{n})|/|(f^{\circ n})'(x_{n})|$ equals
the product of ratios $|f'(y_{n-i})|/|f'(x_{n-i})|$ where $i$
runs from $0$ to $n-1$.   This product can be factored into two
products,
\[ \prod_{x_{i}, y_{i} \in {\cal V}} \frac{|f'(y_{i})|}{|f'(x_{i})|}
\hskip7pt and \hskip7pt
\prod_{x_{i}, y_{i} \in {\cal U}} \frac{|f'(y_{i})|}{|f'(x_{i})|}.\]
We note that the subscript $i$ in
the products are integers in the range $[1, n]$.

\vskip5pt
Using Lemma 3 (the naive
distortion lemma),
we can show that the first product
\[ \prod_{x_{i}, y_{i} \in {\cal V}}
\frac{|f'(y_{i})|}{|f'(x_{i})|}\leq
\exp \Big( \frac{K_{1}}{\beta_{1}}
\sum_{i=0}^{n} |x_{i}- y_{i}|^{\alpha} \Big) .\]

\vskip5pt
The second product
\[ \prod_{x_{i}, y_{i} \in {\cal U}} \frac{|f'(y_{i})|}{|f'(x_{i})|}\]
can be factored into three products
\[ \prod_{x_{i}, y_{i}\in {\cal U}}\frac{|h'(y_{i})|}
{|h'(x_{i})|}\cdot
\prod_{x_{i}, y_{i}\in {\cal U}}\frac{|\tilde{f}'(h(y_{i}))|}
{|\tilde{f}'(h(x_{i}))|}\cdot
\prod_{x_{i},y_{i}\in {\cal U}}\frac{|h'(f(x_{i}))|}
{|h'(f(y_{i}))|},\]
by using the formula
\[ f'(x)=\frac{h'(x)\tilde{f}'(h(x))}{h'(f(x))}.\]

\vskip5pt
By using Lemma 3 again,
the first product
\[ \prod_{x_{i}, y_{i}\in {\cal U}}\frac{|h'(y_{i})|}
{|h'(x_{i})|} \leq
\exp \Big( \frac{K_{3}}{\beta_{3}}
\sum_{i=0}^{n} |x_{i}-y_{i}| \Big) \]
 and the second product
\[ \prod_{x_{i}, y_{i}\in {\cal U}}\frac{|\tilde{f}'(h(y_{i}))|}
{|\tilde{f}'(h(x_{i}))|}\leq
\exp \Big( \frac{K_{3}^{\alpha}K_{2}}{\beta_{2}}
\sum_{i=0}^{n} |x_{i}- y_{i}|^{\alpha } \Big) .\]

\vskip5pt
Suppose
$x_{i},$
$y_{i}$ and $c_{k(i)}$ are in the same set $U_{k(i)}$ and
$v_{k(i)}=f(c_{k(i)})$ is the  critical value.
Because $h'(x)=
1/|x-v_{k(i)}|^{\tau_{k(i)}}$ on a
neighborhood of $v_{k(i)}$, where $\tau_{k(i)}=1-1/\gamma_{k(i)}$ and
$\gamma_{k(i)}$ is the exponent
of $f$ at $c_{k(i)}$, the third product
has the form
\[ \prod_{x_{i}, y_{i}\in {\cal U}}
\Big(
\frac{|y_{i-1}-v_{k(i)}|}{|x_{i-1}-v_{k(i)}|}\Big)^{\tau_{k(i)}}.\]
We note that $x_{i-1}=f(x_{i})$ and $y_{i-1}=f(y_{i})$ are the points
near the
 critical value $v_{k(i)}$ for $x_{i}$ and $y_{i}$ are in the set
$U_{k(i)}$.

\vskip5pt
To control the third product
we write
\[ \frac{|y_{i-1}-v_{k(i)}|}
{|x_{i-1}-v_{k(i)}|} =
|1+ \frac{y_{i-1}-x_{i-1}}
{x_{i-1}-v_{k(i)}}|,\]
which is less than or equal to
$1+|x_{i-1}-y_{i-1}|/|x_{i-1}-v_{k(i)}|$,
for every pair $x_{i}$ and $y_{i}$ in ${\cal U}$.

\vskip5pt
Suppose $l$ is the smallest positive integer such that
$x_{l}$
and $y_{l}$ are in ${\cal U}$.
We consider $l$ in the two cases. The first case is that $l=1$ and
the second case is that $l>1$.

\vskip5pt
In the first case, the images of $x_{l}$ and $y_{l}$ under $f$
are $x$ and $y$.  We have that
\[
\frac{|x-
y|}
{|x-v_{k(l)}|} \leq
\frac{|x- y|}{D_{xy}}.\]

\vskip5pt
In the second case,
suppose $I_{l}$ is the smallest interval containing $x_{l}$, $y_{l}$ and
$c_{k(l)}$ and $I_{l-i}=f^{\circ i}(I_{l})$ for $i=0$, $\cdots $, $l$.
Because the intervals $I_{l-i}$ are contained in ${\cal V}$ for $i=1$,
$\cdots $, $l-1$ (we can always reduce to this case), by using $(VI)$
of Definition 4 and Lemma 3, there is a constant $K_{4}>1$ such that

\[ \frac{ |x_{l-1}-
y_{l-1}|}{|x_{l-1}-v_{k(l)}|}\leq
K_{4} \frac{ |x- y|}{ |x- f^{\circ (l-1)}(v_{k(l)})|}.\]
We note that
$f^{\circ l}(x_{l})=x$ and
$f^{\circ l}(y_{l})=y$ (see Figure 4-a).
This implies that
\[ \frac{ |y_{ l-1 }-x_{l-1}|}
{|x_{ l-1 }-v_{k(l)}| }
\leq
K_{4}
\frac{ |x-y|}
{|x-f^{\circ (l-1)}(v_{k(l)})|}\leq
K_{4}\frac{|x-y|}
{D_{xy}}.\]

\vskip5pt
\psfig{figure=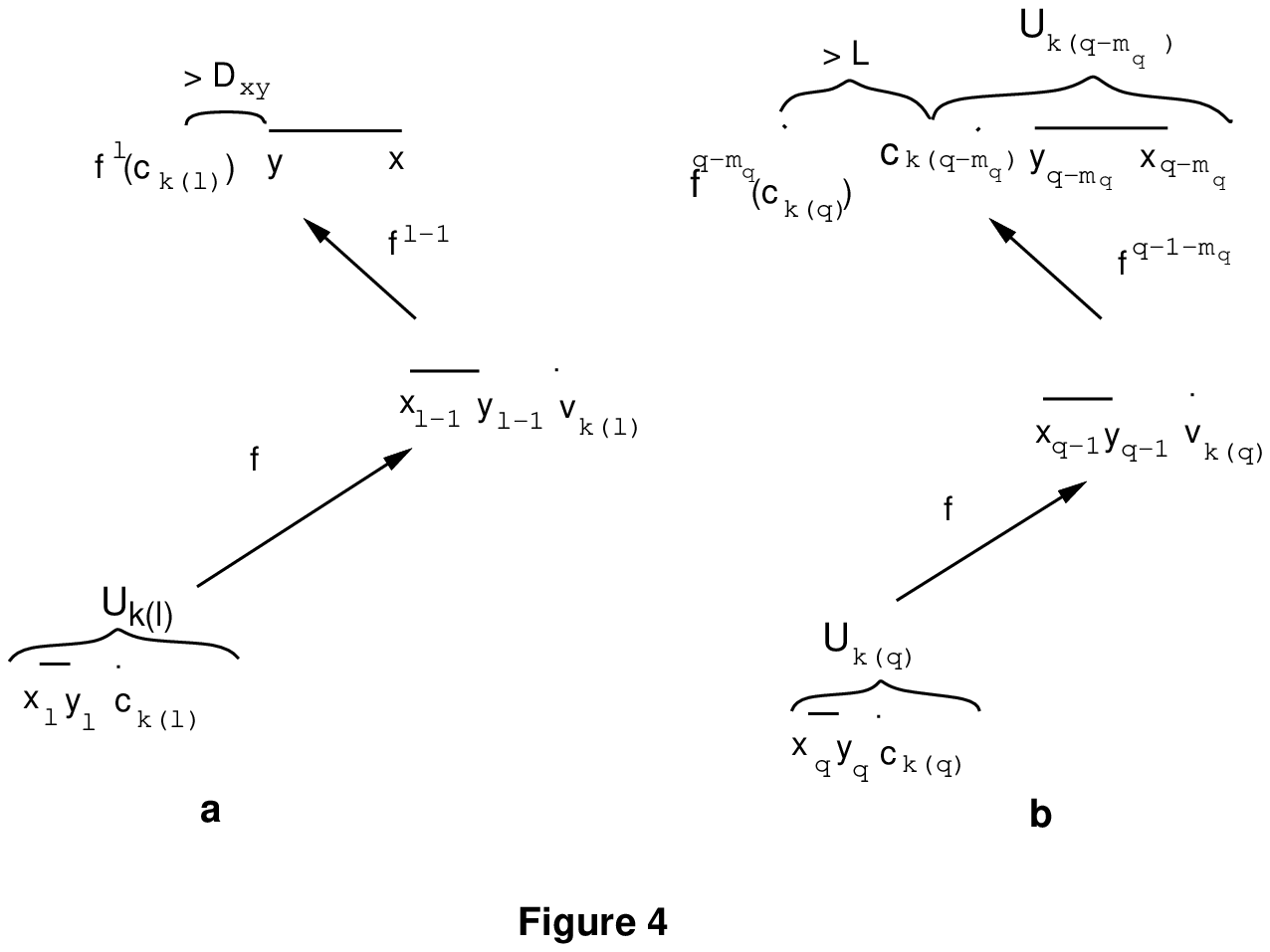}

\vskip5pt
For any $q> l$ with $x_{q}$ and
$y_{q}$ in ${\cal U}$,
let $m_{q}$ be the smallest positive integer such that
$x_{q-m_{q}}$ and $y_{q-m_{q}}$ are in ${\cal U}$
(see Figure 4-b).

\vskip5pt
Suppose $I_{q}$ is the smallest interval containing $x_{q}$, $y_{q}$ and
$c_{k(q)}$ and $I_{q-i}=f^{\circ i}(I_{q})$ for $i=0$, $\cdots $,
$m_{q}$. The intervals $I_{q-i}$ for $i=1$, $\cdots$, $m_{q}-1$ are
contained in ${\cal V}$ (we always can reduce to this case).
By using $(VI)$ of Definition 4 and Lemma 3, there is a positive
constant, we still denote it as $K_{4}$, such that

\[ \frac{
|y_{q-1}-x_{q-1}|}
{|x_{q-1}-v_{k(q)}|}
\leq
K_{4}\frac{ |y_{q-m_{q}}-x_{q-m_{q}}|}{
|x_{q-m_{q}}-
f^{\circ (q-m_{q})}(c_{k(q)})|}.\]
Because $x_{q-m_{q}}$ is in ${\cal U}$ and $f^{\circ
(q-m_{q})}(c_{k(q)})$
is not in ${\cal U}$, the number $|x_{q-m_{q}}-f^{\circ
(q-m_{q})}(c_{k(q)})|$ is bigger than or equal
to $L$, the distance between the set ${\cal U}$
and the closure of the post-critical orbits $\cup_{n=1}^{\infty }
f^{\circ
n}(CP)$. Hence we get

\[ \frac{|x_{q-1}-y_{q-1}|} {|x_{q-1}-v_{k(q)}|}
\leq K_{4}\frac{|x_{q-m_{q}}-y_{q-m_{q}}|}{L}.\]

\vskip5pt
Now the third product satisfies that
\[ \prod_{x_{i}, y_{i}\in {\cal U}}
\Big( \frac{|y_{i-1}-v_{k(i)}|}
{|x_{i-1}-v_{k(i)}|} \Big)^{\tau_{k(i)}}
\leq \exp \Big( \frac{ K_{4}|x-y|}{\tau D_{xy}}
+\frac{K_{4}}{L\tau } \sum_{i=1}^{n} |x_{i}-y_{i}| \Big), \]
where $\tau$ is the maximum of $\tau_{j} =1-1/\gamma_{j}$ for $j=1$,
$\cdots $, $d$.

\vskip5pt
We now prove Lemma 2 by
putting all the estimates together and
$A=
K_{1}/c_{1}+(K_{3}^{\alpha }K_{2})/c_{2}+
K_{3}/c_{3} +
K_{4}/(L\tau)$ and
$B=
K_{4}/\tau$.

\vskip5pt
\noindent {\large \bf \S 3.3
A larger class of one-dimensional mappings.}

\vskip5pt
We can actually prove Lemma 2 for a wider class of one-dimensional
mappings as follows.

\vskip5pt
Suppose $f:M\mapsto M$ is a $C^{1}$-mapping with only
power law critical points.  Let $CP=\{ c_{1}$, $\cdots $, $c_{d}\}$
be the set of critical points of $f$ and $\Gamma=\{
\gamma_{1}$, $\cdots$, $\gamma_{d}\}$ be the set of corresponding
exponents.  Suppose $\eta_{0}$ be the set of the closures of the
intervals of the complement of the set $CP$ of critical points of $f$ in
$M$.

\vskip5pt
{\sc Definition 5.} {\em We say $f$ is $C^{1+\alpha}$ for some $0< \alpha
\leq 1$ if

\vskip5pt
$(1)$ the restriction of $f$ to every interval in $\eta_{0}$ is
differentiable with $\alpha $-H\"older continuous derivative,

\vskip5pt
$(2)$ for every  critical point $c_{i}$, there is a neighborhood
$U_{i}$ of $c_{i}$ such that the functions
$r_{i,-}(x)=f'(x)/|x-c_{i}|^{\gamma_{i}-1}$ for $x<c_{i}$ in $U_{i}$ and
$r_{i,+}(x)=f'(x)/|x-c_{i}|^{\gamma_{i}-1}$ for $x>c_{i}$ in $U_{i}$ are
$\alpha$-H\"older continuous.}

\vskip5pt
Suppose ${\cal U}$ is the union of $U_{i}$ for $i=1$, $\cdots$, $d$
and ${\cal V}$ is the closure of the complement of ${\cal U}$ in $M$.
Let $U_{i-}$ be the subset consisting of all points $x$ in $U_{i}$ with
$x\leq c_{i}$ and
$U_{i+}$ be the subset consisting of all points $x$ in $U_{i}$ with
$x\geq c_{i}$, for
$i=1$, $\cdots $, $d$. Suppose ${\cal W}$ be the collection of all
$U_{i-}$ and $U_{i+}$.

\vskip5pt
{\sc Definition 6.} {\em Suppose $f:M\mapsto M$ is a $C^{1}$-mapping.
We say $f$ is a good $C^{1+\alpha}$-mapping (or good mapping) for some
$0<\alpha \leq 1$ if

\vskip3pt
$(I)$ $f$ has only finitely many power law critical points,

\vskip3pt
$(II)$ $f$ is $C^{1+\alpha}$,

\vskip3pt
$(III)$ there is a positive integer $N$ such that the set $CP$ of
critical points and the closure of the set
$\cup_{n=N}^{\infty}f^{\circ n}(CP)$ are disjoint,

\vskip3pt
$(IV)$ there are two constants $K>0$ and $\nu >1$ such that
for any ${\cal O}_{x, n}=\{ x$, $f(x)$, $\cdots $, $f^{\circ
(n-1)}(x)\}$
with ${\cal O}_{x,n}\cap {\cal U}=\emptyset$,
$|(f^{\circ
k})'(x)|\geq K\nu^{k}$ for any $1\leq k\leq n$.}

\vskip5pt
We say a sequence ${\cal
I}=\{ I_{j}\}_{j=0}^{n}$ of intervals of $M$ is suitable if

\vskip3pt
$(i)$ $I_{j}$ is
the image of $I_{j+1}$ under $f$ for $j=0$, $\cdots$ $n-1$ and

\vskip3pt
$(ii)$ either $I_{j}$ is in ${\cal V}$
or $I_{j}$ is in some interval in ${\cal W}$, for every $j=0$, $\cdots
$, $n$.

\vskip5pt
{\sc Lemma 4} ($C^{1+\alpha}$-Denjoy-Koebe distortion lemma). {\em
Suppose $f$ is a
good $C^{1+\alpha}$-mapping for some $0<\alpha \leq 1$. There are
positive constants $A$ and $B$ such that for
any
suitable sequence ${\cal I}=\{ I_{j}\}_{j=0}^{n}$ of intervals of $M$
and any pair $x$ and $y$ in $I_{0}$,
the distortion of $f$ at $x$ and $y$ along ${\cal I}$
satisfies
 \[
\frac{|g_{n}'(x)|}{|g_{n}'(y)|}\leq \exp
\Big( A\sum_{i=0}^{n}|x_{i}-y_{i}|^{\alpha}
+\frac{B|x-y|}{D_{xy}}\Big) \]
where $D_{xy}$ is the distance between the set $\{x,$ $y\}$ and
the post-critical orbit $\cup_{n=1}^{\infty}f^{\circ n}(CP)$.
}

\vskip5pt
The idea of the proof of this lemma is the same as that of Lemma 2. Details
will be omitted.

\end{document}